\documentstyle[11pt,psfig]{amsart}

\newtheorem{theorem}{Theorem} 
\newtheorem{proposition}[theorem]{Proposition}
\newtheorem{lemma}[theorem]{Lemma}

\newenvironment{proof}{\trivlist\item[\hskip \labelsep {\sc
  Proof:}\enskip]}{\unskip\nobreak\hskip 2em plus 1fil\nobreak%
\hfill\fbox{\rule{0ex}{.5ex}\hspace{.5ex}\rule{0ex}{.5ex}}\endtrivlist}

\newcommand{\R}{{\Bbb R}}
\newcommand{\C}{{\Bbb C}}
\newcommand{\Q}{{\Bbb Q}}
\newcommand{\Z}{{\Bbb Z}}
\renewcommand{\d}{\partial}
\newcommand{\st}{{\bigm|}}
\newcommand{\bq}{\begin{equation}}
\newcommand{\eq}{\end{equation}}
\newcommand{\opname}[1]{\mathop{\fam0#1}} 
\newcommand{\im}{\opname{Im}}
\renewcommand{\hat}[1]{\widehat{#1}}
\renewcommand{\bar}[1]{\overline{#1}}
\renewcommand{\tilde}[1]{\widetilde{#1}}

\newcommand{\PL}{\opname{PL}}
\newcommand{\ie}{{\em i.e.}}
\newcommand{\Ie}{{\em I.e.}}

  \newcommand{\fig}[4]{
    \begin{figure}[ht]
    \centerline{\psfig{file=#1.eps,height=#3in}}
    \caption{#4}\label{#1}
    \end{figure}}
\def\fcorner {\fig{}{2.21}{2}
    {Corner separation}}
\def\fmirror {\fig{}{4.04}{1.8}
    {Mirror image}}
\def\finsert {\fig{}{4.19}{2.28}
    {Insertion}}
\def\fbridge {\fig{}{4.03}{1.65}
    {Bridge immersion}}
\def\fpartial {\fig{}{2.04}{2.04}
    {Partial insertion}}
\def\fentryexit {\fig{}{4.78}{1.86}
    {Self-insertion: leaves a \& b}}
\def\frainbow {\fig{}{3.65}{.5}
    {Matched transitions}}
\def\fwilson {\fig{}{3.1}{2.45}
    {The plug ${{\cal W}}$}}
\def\fdad {\fig{}{4}{2.7}
    {An analytic plug}}
\def\ftoungue {\fig{}{2.71}{2.82}
    {Bridge immersion of $I\times T^2$}}
\def\fcollar {\fig{}{3.28}{1.78}
    {PL vertical collar}}
\def\fvertical {\fig{}{4.31}{2.15}
    {A 2-dimensional minimal set}}
\def\fmapf {\fig{}{3.28}{1.40}
    {The map $f$}}
\def\fsigma {\fig{}{5.07}{2.07}
    {$\sigma $ restricted to $A$}}
\def\f1dim {\fig{}{2.17}{2.30}
    {A 1-dimensional minimal set}}

\def\picture #1 by #2 (#3){
  \vbox to #2{
    \hrule width #1 height 0pt depth 0pt
    \vfill
    \special{picture #3} % this is the low-level interface
}
}  

\author{Greg Kuperberg}
\address[G. Kuperberg]{Department of Mathematics \\
        University of California, Davis \\ Davis, CA 95616}
\email{greg@@math.uchicago.edu}
\thanks{The first author was supported by an NSF
        Postdoctoral Fellowship, grant \#DMS-9107908.}
\author{Krystyna Kuperberg}
\address[K. Kuperberg]{Department of Mathematics \\ Auburn University \\
        Auburn, AL 36830}
\email{kuperkm@@mail.auburn.edu}
\thanks{The second author was supported in part by NSF grant \#DMS-9401408.}
\title{Generalized counterexamples to the Seifert conjecture}
\date{July 6, 1994} 
\begin{document}

\begin{abstract}
Using the theory of plugs and the self-insertion construction due
to the second author, we prove that a foliation of any codimension of
any manifold can be modified in a real analytic or piecewise-linear fashion
so that all minimal sets have codimension 1.  In particular,
the 3-sphere $S^3$ has a real analytic dynamical system such that
all limit sets are 2-dimensional. We also prove that a 1-dimensional
foliation of a manifold of dimension at least 3 can be modified in
a piecewise-linear fashion so that so that there are no closed leaves but
all minimal sets are 1-dimensional.  These theorems provide new
counterexamples to the Seifert conjecture, which asserts that every
dynamical system on $S^3$ with no singular points has a periodic
trajectory.
\end{abstract}

\maketitle

\section{Introduction}

In 1950, H.~Seifert \cite{Seifert} asked whether every dynamical system on the
3-sphere with no singular points has a periodic trajectory.  The conjecture that
the answer is yes became known as the Seifert conjecture.  
Seifert proved the conjecture for perturbations of the flow parallel
to the Hopf fibration.
In 1974,
P.~A.~Schweitzer \cite{Schweitzer} found a $C^1$ counterexample to the Seifert
conjecture, which was modified to a $C^2$ counterexample by J.~Harrison
\cite{Harrison}.
In 1993, H.~Hofer \cite{Hofer} proved the conjecture for contact flows.

The main idea of this paper comes from the construction of a smooth
counterexample to the Seifert conjecture due to the second author
\cite{kseifert}.   Here, we outline more general consequences of that
idea.  In particular, we establish the following theorem:

\begin{theorem}  A foliation of any codimension of any manifold
can be modified in an analytic or piecewise-linear fashion so
that all minimal sets have codimension 1.
\label{thcd1}
\end{theorem}

This theorem strengthens an earlier result of F.~W.~Wilson
\cite{Wilson} which for oriented 1-foliations, establishes minimal sets of
codimension 2.

Theorem~\ref{thcd1} is trivial for foliations which are themselves codimension
1.  The problem of opening all compact leaves in codimension 1, which is a more
natural question, has been partly settled by S.~P.~Novikov\cite{Novikov} and
P.~A.~Schweitzer \cite{Paul}. Novikov proved that every $C^2$ codimension
1 foliation of $S^3$ has a closed leaf (which was later extended to continuous
foliations by V.~V.~Solodov \cite{Solodov} and  G.~Hector and U.~Hirsch
\cite{HH}), while Schweitzer has shown that it is possible to modify any
codimension 1 foliation in dimension 4 or higher in a $C^1$ fashion so that
it has no compact leaf.

Here and throughout the paper, smooth means $C^\infty$ and
analytic means real analytic or $C^\omega$.  All manifolds are assumed to be
paracompact and Hausdorff, and they are assumed to have no boundary
unless explicitly stated otherwise.  Also, dimension means
covering dimension, cohomological dimension, small inductive
dimension, or large inductive dimension; they are all equal for
topological spaces considered in this paper.
In our context, we can freely convert between an (oriented)
1-dimensional foliation and a dynamical system with no singular points. 
Note that
Theorem~\ref{thcd1} implies the following theorem:

\begin{theorem} The 3-sphere $S^3$ has an analytic dynamical
system such that all limit sets are 2-dimensional.  In
particular, it has no circular trajectories.  \label{th3sph} 
\end{theorem}

Table~\ref{texamples} shows the best known continuity of various
kinds of foliations of 3-manifolds.  Entries with no citation
are covered by Theorem~\ref{th3sph} or its analogue for arbitrary
3-manifolds.

\begin{table}[ht]
\begin{center}
\begin{tabular}{c||c|c}
& not volume-preserving & volume-preserving \\ \hline \hline
discrete circles   & $C^\omega \cite{Wilson}$ & 
$C^\infty$ \& $\PL$ \cite{kvol} \\ \hline
1-dim. minimal sets & $C^1$ \cite{Schweitzer}, $C^2$ \cite{Harrison}, $\PL$ &
$C^1$ \cite{kvol} \\ 
but no circles & & \\\hline
2-dim. minimal sets & $C^\omega$, $\PL$ & --- 
\end{tabular}
\vspace{\baselineskip}
\caption{Known foliations of 3-manifolds}
\label{texamples}
\end{center}
\end{table}

All counterexamples to the Seifert conjecture and its analogues described in
this paper are based on constructions of aperiodic plugs.  An (insertible, untwisted,
attachable) plug, whose prototype was defined by Wilson, is an oriented, 1-dimensional
foliation ${\cal F}$ of a manifold with boundary, the Cartesian product $F \times I$ of an
$(n-1)$-dimensional manifold $F$ and $I$.   The foliation ${\cal F}$ agrees with
the trivial foliation in the $I$ direction on a neighborhood of $\d (F \times
I)$, if a leaf connects $(p,0)$ with $(q,1)$, then $p = q$, and there is a
non-compact leaf containing some $(p,0)$. Here the base $F$ is an
$(n-1)$-manifold that admits a bridge immersion in $\R^{n-1}$, an immersion that
lifts to an embedding $\R^n$.   A plug is aperiodic if it has no closed leaves. 
For example, the 3-dimensional Wilson plug on $A \times I$, where $A$ is an
annulus, is not aperiodic,  while the Schweitzer plug is an aperiodic plug on
$pT \times I$, where $pT$ is a punctured torus. Although Wilson suggested the
technique of inserting plugs to modify foliations, it was Schweitzer who first
used plugs  to break circular leaves, and it was his important observation that
the base of a plug need only admit a bridge immersion, rather than an embedding,
in order to be insertible. The results in this paper are based on an idea for
constructing aperiodic plugs from a Wilson-type plug which breaks its own
circular leaves by means of self-insertion.

The Schweitzer-Harrison $C^2$ aperiodic plug has two 1-dimensional minimal sets. In
the $\PL$ category, we obtain an aperiodic plug with one 1-dimensional minimal set.

\begin{theorem} A 1-foliation  of a manifold of dimension at least 3 can be
modified in a $\PL$ fashion so there are no closed leaves but all minimal sets are
1-dimensional. Moreover, if the manifold is closed, then there is an aperiodic $\PL$
modification with only one minimal set, and the minimal set is 1-dimensional.
\label{th4sph}
\end{theorem}

The authors would like to thank William P. Thurston for observing that the basic
construction is analytic and not merely smooth.

\section{Preliminaries}

This paper will consider four kinds of functions between manifolds: 
continuous,  smooth (meaning $C^\infty$), analytic (meaning real analytic or
$C^\omega$), and piecewise linear or PL.  Each of these four classes of
functions is a {\em smoothness category\/}.

Recall that many kinds of manifolds can be understood by gluing charts. A
smooth manifold has smooth gluing maps, an analytic manifold has analytic
gluing maps, and in general a manifold is said to be in a given smoothness
category if its gluing maps are.  We will need the following fundamental
result \cite{Morrey} \cite{Grauert}:

\begin{theorem}[Morrey, Grauert] Two analytic manifolds which are
diffeomorphic  are analytically diffeomorphic.
\end{theorem}

The Morrey-Grauert theorem has many formulations, all of which were rendered
equivalent by Whitney \cite{Whitney} in work that predated the theorem
itself: Given points $p$ and $q$ in $M$, there exists an analytic function
$f$ such that $f(p) \ne f(q)$. (This formulation is the closest to what
Morrey and Grauert proved directly.) Every analytic manifold $M$ admits an
embedding in some $\R^n$. Given a $C^n$ function $f$ on $M$, there exists a
sequence $\{f_i\}$ of analytic functions on $M$ such that $f_i$ and its
first $n$ derivatives converge pointwise to those of $f$.  Another relevant
result due to Whitney is the fact that every smooth manifold admits a
smoothly compatible real analytic structure.

A {\em $k$-dimensional foliation structure\/} or {\em $k$-foliation\/} on an
$n$-manifold $M$ is an atlas of charts in $\R^n$ that preserve the parallel
$k$-plane foliation of $\R^n$, which is a partition of $\R^n$ into translates
of flat $\R^k \subset \R^n$.  $M$ is then a {\em $k$-foliated manifold\/}.
The foliation structure is in a given category, such as smooth, if the gluing
maps are simultaneously in the same category and preserve $k$-planes.  An
{\em oriented foliation\/} is a foliation structure whose gluing maps
preserve the standard orientation of $\R^k$ and its translates; however, they
need not preserve the orientation of $\R^n$.

Let $M$ be a $k$-foliated manifold.  Consider the topology on $\R^n = \R^k
\times \R^{k-n}$ which is the usual topology in the direction of the first
factor $\R^k$ and the discrete topology in the direction of the second factor
$\R^{n-k}$.  This topology can be restricted to a topology on charts in $\R^n$
and then pushed forward to a topology on $M$ using the charts that describe the
foliation of $M$.  The resulting topology is the {\em leaf topology\/} on $M$
and it divides $M$ into a disjoint union of connected $k$-manifolds which are
the {\em leaves\/} of the foliation of $M$.  A foliation is often described in
terms of its leaves.  Also, a map between foliated manifolds is {\em
leaf-preserving\/} if it is continuous in both the leaf topology and the usual
topology. (Unless explicitly stated otherwise, all terms such as closed,
connected, etc., will refer to the usual topology on $M$ rather than the leaf
topology.)  A {\em minimal set\/} of a foliation is a set which is minimal among
non-empty, compact subsets which are unions of leaves.  A smooth or analytic
foliation is also uniquely determined by the $k$-plane field ($k$-dimensional
subbundle of the tangent bundle) parallel to it.

If $M$ is $n$-manifold with a $k$-foliation ${\cal F}$ and a $j$-submanifold $N$,
then $N$ and ${\cal F}$ are {\em transverse\/} at a point $p$ if there is a local
equivalence that sends a neighborhood $U$ of $p$ into $\R^n$, $N\cap U$ into a
$j$-plane, and ${{\cal F}}|_U$ into parallel $k$-planes that intersect the
$j$-plane in $(j+k-n)$-planes.  In particular, if $N$ is transverse to a
foliation, then $N$ admits a tubular neighborhood; transversality is problematic
without this condition.  Note that transversality is defined relative to its
smoothness category; for example, in the plane the curves $y = x^3$ and $y = 0$
are continuously but not smoothly transverse.

A {\em flow\/} on a metric space $A$ is a system of partial homeomorphisms of
$A$ parametrized by time $t \in \R$.  Specifically, let $U$ be an open subset
 of $\R \times A$ containing $\{0\} \times A$ such that $U \cap (\R \times
\{p\})$ is connected for every $p$.  A flow on $A$ is a continuous map $\Phi:U
\to A$ such that $\Phi(0,p)=p$ and whenever $\Phi(t,p)$ and $\Phi(s,\Phi(t,p))$
are both defined, $\Phi(t+s,p)$ is also defined and equals $\Phi(s,\Phi(t,p))$.
If $U = \R \times A$, then $\Phi$ is a {\em dynamical system\/}, or equivalently
a {\em topological group action\/} of $\R$ on $A$.  For example, a flow on a
closed manifold is necessarily a dynamical system. If $\Phi$ is either a flow or
a dynamical system, a {\em trajectory\/} of a point $p$ is the image of
$\Phi((\R \times \{p\}) \cap U)$ in $A$. A point $p$ is a {\em rest point\/} if
its trajectory consists only of $p$. A {\em positive\/} ({\em negative\/}) {\em
limit set\/} of a trajectory of $p$ is the set of limit points of sequences
$\{\Phi (t_n,p)\},$ such that $t_n\rightarrow \infty$ ($t_n\rightarrow -\infty)$.

If $M$ is a manifold and $\Phi$ is smooth or analytic, then  ${d\Phi \over
dt}|_{t = 0}$ is the vector field of $\Phi$ and is in the same smoothness
category as $\Phi$. Contrariwise, a standard integration theorem for
differential equations says that a $C^1$ vector field can be integrated to
produce a corresponding flow or dynamical system. If $\Phi$ has no rest points,
the collection of all of its trajectories is the set of leaves of an oriented
1-dimensional foliation which is in the same smoothness category as $\Phi$. 
Contrariwise, if ${\cal F}$ is an oriented, 1-dimensional foliation which is
smooth or analytic, there exists a parallel non-vanishing vector field in this
same category as ${\cal F}$ which can be integrated to produce a compatible flow. 
For example, a smooth or analytic manifold always admits a Riemannian metric in
the same category, and the unit vector field parallel to the foliation and
pointing in the direction of its orientation suffices.

An oriented 1-foliation ${\cal F}$ which is PL or continuous  always has a
parallel flow, by the following constructions:  In the PL case, embed the
manifold $M$ of ${\cal F}$ in some $\R^n$, and define $\Phi(t,x)$ so that the
oriented length of the leaf segment from $x$ to $\Phi(t,x)$ is $t$.  For the
continuous case choose a locally finite atlas of charts $\alpha_i:\R^n \to M$,
and choose a partition of unity $\{f_i\}$ that refines $\{\alpha_i\}$.  Let $s$
be a line segment in some leaf with endpoints $x$ and $y$ which is oriented from
$x$ to $y$, and suppose that $s$ is entirely contained in every chart which
intersects it.  Then define $t$ by the formula $$t = \sum_i
\int_{\alpha_i^{-1}(s)} f_i$$ and define $y = \Phi(t,x)$.  The function $\Phi$
then extends uniquely to a flow with no rest points.

Henceforth, we will loosely switch between flows, vector fields, and oriented
1-foliations, since they are almost the same geometrically.

If $M$ is a manifold with boundary, a foliation  of $M$ is a foliation of the
interior of $M$ which admits an extension to a foliation of an open manifold
containing $M$; similarly, a flow on $M$ is a flow on the interior of $M$ which
admits an extension to a flow on an open manifold containing $M$.  If ${\cal F}$ is
a  1-foliation, then the {\em parallel boundary\/} of ${\cal F}$ is the subset of
$\d M$ where ${\cal F}$ is locally modelled by the foliation of upper half space by
horizontal lines; similarly the {\em transverse boundary\/} is the subset of $\d
M$ where ${\cal F}$ is locally modelled by the foliation by vertical lines.  Note
that $M$ can have boundary which is neither parallel nor transverse. A  {\em
minimal set  of a flow\/} $\Phi$ on a manifold with boundary is  a set $A$ such
that  the flow restricted to $A$ is a dynamical system, and $A$ is a minimal set
of this dynamical system. A leaf  may have boundary or {\em endpoint(s)\/}. An
{\em infinite leaf\/}  is a non-compact leaf. A leaf that is neither closed nor
infinite is a {\em finite leaf\/}.

\fcorner

A {\em manifold with corners} is a smooth or analytic manifold with piecewise
smooth or analytic boundary.  \Ie, let $M$ be a compact $n$-manifold with
boundary. If, for each point $p \in M$, there exists an $n$-dimensional PL
submanifold $P$ of $\R^n$ with boundary which is diffeomorphic to a neighborhood
of $p$.  For example, a manifold with boundary is a manifold with corners, and
the Cartesian product of two manifolds with corners is a manifold with corners.
Flows and foliations on manifolds with corners are defined in the same way as on
manifolds with boundary, and the standard boundary types of parallel and
transverse boundary can be extended to manifolds with corners.  Let $N$ be a PL
$(n-1)$-submanifold of $\R^{n-1}$.  A 1-foliation ${\cal F}$ has {\em corner
separation\/} at $p \in \d M$  if $p$ is in neither parallel nor transverse
boundary, but if ${\cal F}$ is locally equivalent at $p$ to $N \times [0,\infty)$
foliated by rays $\{x\} \times [0,\infty)$. In particular, $p$ belongs to the
closure of the transverse boundary and the closure of the parallel
boundary. Figure~\ref{fcorner} gives an example of corner separation between
parallel and transverse boundary.

Although the Morrey-Grauert theorem does not directly apply to a
manifold with corners $M$, a manifold with corners is always contained
in an open manifold, and the Morrey-Grauert theorem applied to this
larger manifold implies that analytic functions on $M$ separate points
and that $M$ admits an analytic embedding with the same geometry
as any smooth embedding.

\section{\label{splugs} Plugs}

Except where explicitly stated otherwise, the constructions in this section apply
uniformly in each of the four smoothness categories.  Typically, an object $O$
might involve a foliation, a gluing map, and an embedding, in which case $O$ is
a smooth object (for example) if all three parts are in the smooth category.

A {\em flow bordism} is an oriented 1-foliation ${\cal P}$ of a connected, compact
manifold $P$ with boundary or corners such that $\d P$ is entirely transverse
boundary, parallel boundary, or corner separation, and such that all leaves in
the parallel boundary of $P$ are line segments. (A flow bordism is very similar
to an oriented {\em foliated cobordism} \cite{Tamura}.  The principal differences
are that a foliated cobordism may be a higher-dimensional foliation and it is
usually considered without corner separation or parallel boundary.)  If ${\cal P}$
is a flow bordism, let $F_-$ be the closure of the transverse boundary oriented
inward, and similarly let $F_+$ be the closure of the transverse boundary
oriented outward.  The foliation ${\cal P}$ might in addition have one or both of
the following properties:

\begin{description} 
\item[(i)] There exists
an infinite leaf with an endpoint in $F_-$. 
\item[(ii)] There exists a
manifold $F$ and homeomorphisms $\alpha_\pm:F \to F_\pm$ such that if
$\alpha_+(p)$ and $\alpha_-(q)$ are endpoints of a leaf of ${\cal P}$,
then $p=q$. 
\end{description}

If ${\cal P}$ satisfies property {\bf (ii)}, it has {\em matched ends\/}. The foliation
${\cal P}$ is a {\em plug\/} if it has properties {\bf (i)} and {\bf (ii)}, but
only a {\em semi-plug\/} if it has property {\bf (i)} but not property {\bf
(ii)}.  It is an {\em un-plug\/} if has property {\bf (ii)} but not property
{\bf (i)}.  Table~\ref{tplugtypes} gives a summary of these four definitions.
The manifold $F_-$ is the entry region of ${\cal P}$, while $F_+$ is the exit
region. If ${\cal P}$ has matched ends, then $F$ is the {\em base\/} of $\cal P$. 
The manifold $P$ is the {\em support} of ${\cal P}$.  The {\em entry stopped
set\/} $S_-$ of ${\cal P}$ is the set of points of $F_-$ which  are endpoints of
infinite leaves; the {\em exit stopped set\/} $S_+$ is defined similarly.  If
${\cal P}$ has  matched ends, the stopped set $S$ is defined as
$\alpha_-^{-1}(S_-) = \alpha_+^{-1}(S_+)$. If $S$ has non-empty interior, then 
${\cal P}$ {\em stops content}.

\begin{table}[ht]
\begin{center}
\begin{tabular}{c||c|c|c}
          &  not {\bf (ii)}     &  {\bf (ii)}    \\ \hline \hline
not {\bf (i)}  & semi-un-plug & un-plug \\ \hline
{\bf (i)}      & semi-plug    & plug
\end{tabular}
\vspace{\baselineskip}
\caption{Describing ${\cal P}$ based on its properties}
\label{tplugtypes}
\end{center}
\end{table}

Note that if ${\cal P}$ is a flow bordism with entry region $F_-$ and exit region
$F_+$, then since $F_-$ and $F_+$ are  the transverse boundaries together  with
some corner separation points, they possess trivially foliated neighborhoods. 
Specifically, if $F_\pm \times [0,1)$ is foliated by vertical fibers $\{p\}
\times [0,1)$, there exists a leaf-preserving homeomorphism $\omega_\pm$ from an
open neighborhood of $F_\pm$ to $F_\pm \times [0,1)$; moreover, we can take
$\omega_\pm(p) = (p,0)$ for all $p \in F_\pm$.

\begin{lemma} If ${\cal P}$ is a flow bordism with connected support $P$ and $\cal
P$ has at least one circular or infinite leaf, then the entry stopped set $S_-$
(resp. the exit stopped set $S_+$) is non-empty if and only
if $F_-$ (resp. $F_+$) is non-empty.  In
particular, if ${\cal P}$ has matched ends (and $F_- \cong F_+ \cong F$ is
non-empty), then ${\cal P}$ is a plug. \label{lisplug}  
\end{lemma} 
\begin{proof}
Let $\Phi$ be a flow parallel to ${\cal P}$.  If $p \in P$, then $L_+(p)$, the
{\em future longevity\/} of $p$, is the supremum of $t$ such that $\Phi(t,p)$ is
defined.  Similarly, the {\em past longevity\/} $L_-(p)$ is the supremum of $t$ such that
$\Phi(-t,p)$ is defined. Thus we have two functions
$L_\pm: P \to [0,\infty]$. If $L_+ (p) <\infty $
or $L_- (p) <\infty $, then the leaf containing $p$ has an endpoint $q$ on
$F_+$ or $F_-$, respectively. Investigating the points in a neighborhood of
such a
$q$, we see that two sets $L_\pm ^{-1} (\infty ) $ are both closed.  If
$S_-=\emptyset$, then $F_-\cap L_+^{-1}([0,\infty )) = F_-$,
and the set $A$ consisting of leaves that intersect $F_-$ is closed.  Then
$P$ is the union of two disjoint non-empty closed sets, $A$
and $L_- ^{-1} (\infty ) \cup L_+ ^{-1} (\infty ) $,
which is a contradiction.  Hence  $S_-\not=\emptyset$ if
$F_-\not=\emptyset$.
Similarly,  $S_+\not=\emptyset$ if $F_+\not=\emptyset$.
\end{proof}

\fmirror

\begin{lemma} If ${\cal P}$ is a semi-un-plug with support $P$
and entry region $F_-$, then there is a foliation isomorphism 
$\gamma :F_- \times I \to P$, where $F_- \times I$ is foliated
by fibers $\{p\} \times I$.  In particular, if ${\cal P}$ is an un-plug with base $F$,
then there is a foliation isomorphism $\alpha:F \times I \to P$ which
extends the maps $\alpha_{\pm}$.~\label{listrivial}
\end{lemma}
\begin{proof}
Let $\Phi$ be a flow parallel to ${\cal P}$.  Let $L_+$ be the future
longevity as in the proof of Lemma~\ref{lisplug}.  In each smoothness
category other than the PL category, define $\gamma(p,t) = \Phi(p,t/L_+(p))$.
The map $\gamma$ has all of the desired properties.

Unfortunately, the PL category is not closed under the arithmetic operation of
division, so a more complicated argument is necessary.  Let $U$ be the region
in $F_- \times \R$ bounded by $F_- \times \{0\}$ and the graph of $L_+$
restricted to $F_-$.  The map $(p,t) \mapsto \Phi(p,t)$ is a foliation
isomorphism from $U$ to $P$; it remains to find a foliation isomorphism
$\beta$ from $F_- \times I$ to $U$.  To construct $\beta$, let $T$ be a
triangulation with vertex set $V$ of $F_-$ such that $L_+$ is linear on each
simplex of $T$ and such that $T$ is the barycentric subdivision of another
triangulation.  Since $T$ is a barycentric subdivision, there is a map $o:V
\to \{1,\ldots,n+1\}$, where $n$ is the dimension of $F_-$, which is
bijective when restricted to the vertices of any single $n$-simplex of $T$.  If
$K$ is such a simplex, consider a triangulation $\{K_i\}$ of $K \times I$,
such that the simplex $K_i$ has vertices $o^{-1}(\{1,\ldots,i\}) \times
\{0\}$ and $o^{-1}(\{i,\ldots,n+1\}) \times \{1\}$.  Then the map on vertices
of $K \times I$ given by $(p,0) \mapsto (p,0)$ and $(p,1) \mapsto (L_+(p),1)$
extends to a PL foliation isomorphism from $K \times I$ to $(K \times \R)
\cap U$. These maps fit together for different simplices $K$ of $T$, yielding
the desired map $\beta$.
\end{proof}

The following construction due to Wilson \cite{Wilson}  turns  a semi-plug into
a plug.  If ${\cal P}_1$ and $\cal P_2$ are two flow bordisms such that the exit
region of ${\cal P}_1$ is the same as the entry region of $\cal P_2$, their {\em
concatenation} is a flow bordism obtained by identifying trivially foliated
neighborhoods of this shared region. The {\em mirror image} $\bar{{\cal P}}$ of
${\cal P}$ is given by reversing the orientation of the leaves of $\cal P$, which
has the effect of switching the entry and exit regions.  The {\em mirror-image
construction} is the concatenation of ${\cal P}$ with $\bar{\cal P}$; it is easy
to see that the result of this concatenation has matched ends.
(Figure~\ref{fmirror} shows a schematic picture of a mirror-image construction.)

The primary purpose of plugs is the operation of {\em insertion}.  The geometric idea
of insertion is illustrated in Figure~\ref{finsert}. Given a foliation ${\cal X}$
of an $n$-manifold $X$ and an $n$-dimensional plug ${\cal P}$ with base $F$, we
wish to find a leaf-preserving embedding of the trivially foliated $F \times I$
in $X$ and replace it with ${\cal P}$.  In full generality, this procedure raises
three technical issues:  insertibility, attachability, and twistedness.

\finsert

An {\em insertion map} for a plug ${\cal P}$ into a foliation $\cal X$ is an
embedding $\sigma : F \to X$ of the base of ${\cal P}$ which is transverse to $\cal X$. 
Such an insertion map can be extended to an embedding $\sigma: F \times I \to
X$ which takes the fiber foliation of $F \times I$ to ${\cal X}$. (It is
convenient to denote the maps $F \to X$ and $ F \times I \to X$ by the same
$\sigma$ and this should not lead to a missunderstanding.)  An
$n$-dimensional plug ${\cal P}$ is {\em insertible\/} if $F$ admits an
embedding in $\R^n$ which is transverse to vertical lines.  Such an
embedding is equivalent to a {\em bridge immersion\/} of $F$ in $\R^{n-1}$, \ie, an
immersion which lifts to an embedding of $F \times I$ in $\R^n$. For example, if $F$ is
2-dimensional, orientable, and has non-empty boundary, then $F$ necessarily
admits a bridge immersion. Figure~\ref{fbridge} shows a bridge immersion of a
punctured torus $pT$ and a bridge immersion of a surface of genus three with
three punctures. The corresponding embedding of $pT \times I$ is  one of the
main steps in Schweitzer's counterexample to the Seifert
conjecture \cite{Schweitzer}.

\fbridge

A plug is {\em attachable} if every leaf in the parallel boundary is finite.
In general, if $M$ is a manifold with boundary, let $N_M$ denote
an open neighborhood of $\d M$ in $M$.
The next
step in plug insertion is to remove $\sigma((F \times I)-N_{F \times I})$
from $X$ and glue the open lip $\sigma(N_{F \times I})$ to the support $P$ of
an attachable plug ${\cal P}$ by a leaf-preserving homeomorphism $\alpha:N_{F
\times I} \to N_P$, where $N_P$ is a neighborhood of $\d P$.  Moreover, the
identification $\alpha$ should satisfy $\alpha(p,0) = \alpha_-(p)$ and
$\alpha(p,1) = \alpha_+(p)$.  A map $\alpha$ with these properties is an
attaching map for ${\cal P}$.  Let $G$ be the parallel boundary of $P$, so that
$\d P = G    \cup F_- \cup F_+$.  Recall that all of the leaves in $G$  have
two endpoints; this must also be true in a neighborhood $N_G$ of $G$.  It
follows by Lemma~\ref{listrivial} that there exists a leaf-preserving
homeomorphism $N_G \to N_F \times I$.  This equivalence, together
with the matched ends condition,
ensures the existence of an attaching map for any attachable plug or un-plug.

Let $\sigma$ be an insertion map of a plug ${\cal P}$ into a foliation $\cal X$ on
a manifold $X$, and let $\hat{{\cal X}}$ be the foliation on the manifold
$\hat{X}$ resulting from inserting ${\cal P}$ into $\cal X$.  The plug $\cal P$ is
{\em untwisted\/} if the attaching map $\alpha$ extends to a homeomorphism $F \times I
\to P$.  The manifolds $X$ and $\hat{X}$ need not be homeomorphic, but they are
if ${\cal P}$ is untwisted.  If $\cal P$ is a smooth or PL plug, then it is assumed
that the extension of $\alpha$ is smooth or PL just as $\alpha$ is. However, if
${\cal P}$ is analytic, then an analytic extension of $\alpha$ is usually not
possible.  In this case, ${\cal P}$ is untwisted if $\alpha$ admits a smooth
extension.  The resulting manifolds $X$ and $\hat{X}$ have analytic structures
which are a priori only smoothly diffeomorphic; however, the Morrey-Grauert
theorem ensures that they are in fact analytically equivalent.  All plugs in
this paper are assumed to be untwisted unless explicitly stated otherwise.

A plug or a semi-plug is {\em aperiodic} if it has no circular leaves.

One of the main reasons to insert an aperiodic plug into another foliation is to
break a circular leaf.  Let ${\cal P}$ be a plug with stopped set $S$ and base
$F$, and let $\sigma$ be an insertion map into a foliation ${\cal X}$. If
$\sigma(S)$ intersects a circular leaf $l$ of ${\cal X}$, then the remnant of $l$
in $\hat{{\cal X}}$ is an infinite leaf. On the other hand, if $l$ is disjoint
from $\sigma(S)$ but intersects $\sigma(F)$, then there is a leaf $\hat{l}$ in
$\hat{{\cal X}}$ corresponding to $l$ with the same topology as $l$, although the
geometry of $\hat{l}$ can differ from that of $l$.

The simplest interesting target for an insertion is an irrational foliation of
the 3-sphere $S^3$.  Given an irrational number $r \in \R$, the irrational
foliation of slope $r$ is parallel to the vector field $\vec{V}$ on $\{(x,y) \in
\C^2 \st |x|^2 + |y|^2 = 1\}$ given by
$$\vec{V}(x,y) =ix+iry\, .$$ 
The foliation has two circles, and the closure of any other leaf is a torus
lying between these two circles.  Two copies of any aperiodic plug can be
inserted in small regions meeting these circles and breaking them 
\cite{Schweitzer}.  Thus, to find a counterexample to the Seifert conjecture, it
suffices to construct an aperiodic plug.  If the plug is smooth or analytic, the
corresponding foliation of $S^3$ is respectively smooth or analytic.

\section{\label{sstoppage} Global stoppage}

In 1935,  K.~Borsuk \cite{Borsuk} gave an example of a fixed point free
homeomorphism of an acyclic compact subset of $\R ^3$; a solid cylinder with
two narrowing and spiraling tunnels drilled out. His example can be easily
modified to obtain a semi-plug on $D^2\times I$ ($D^n$ is an $n$-cell) with
two circular leaves whose entry and exit stopped sets $S_-$ and $S_+$ are
disks.  A similar construction on a solid torus (Borsuk's example can be
obtained by cutting or unwrapping the torus) was given in 1952 by
F.~B.~Fuller \cite{Fuller}. The consequence of Fuller's result that
$D^{n-1}\times S^1$, $n\geq 3$, admits a dynamical system whose only minimal
set is an $(n-2)$-torus $S^1\times \cdots \times S^1$ contained in a
$D^{n-1}\times \{p\}$ section motivated Wilson's introduction of plugs, which
he used to arrest flows globally.
 
\begin{theorem}[Wilson] Every smooth $n$-manifold of Euler characteristic zero or
non-compact has a smooth dynamical system with a discrete set of minimal sets. Each of
the minimal sets is an $(n-2)$-torus, and every trajectory originates (resp. limits)
on one of these tori.  \label{thwilson}
\end{theorem}

Wilson's theorem implies that an analogue to the Seifert conjecture for higher
dimensional spheres of odd dimension is false.  Although Wilson claimed
Theorem~\ref{thwilson} in the smooth category, the construction is actually
analytic. Also, with slight modifications, Wilson's construction is valid in
a PL setting.

Wilson constructed a content-stopping, mirror-image plug and inserted many small
copies of the plug into the foliation to break the leaves into small pieces. 
More specifically, he proved that if ${\cal F}$ is an oriented 1-foliation of an
$n$-manifold $M$, $n\geq 3$, and ${\cal U}$ is open cover of $M$, then $\cal X$
can be modified to an oriented 1-foliation  with every leaf contained in a star
of an element of ${\cal U}$.  The following result of \cite{Reed} is
therefore a special case of Theorem~\ref{thwilson}.

\begin{theorem} For every $\epsilon > 0$, there exists an oriented 1-foliation
of $R^3$ with all leaves of diameter less than $\epsilon$. \label{threed}
\end{theorem}

To establish an aperiodic version of Theorem~\ref{thwilson} and \ref{threed}
in dimension 3, it suffices to use an aperiodic plug that stops content, for
instance, a variation of Schweitzer's plug \cite{Coke} or a plug  obtained by
breaking the circles of Wilson's periodic plug with any aperiodic plug.
Sections~\ref{scounter}, \ref{shigher}, and \ref{spl} give constructions of
analytic and PL aperiodic plugs such that each has a unique minimal set and
it has codimension 1.  These plugs yield stronger versions of
Theorems~\ref{thwilson} and \ref{threed}.

\begin{theorem} If $M$ is a continuous, $C^r$, $C^\infty$, $C^\omega$, or $\PL$ manifold of
dimension $\geq 3$ admitting an oriented 1-foliation in the same smoothness
category,  and \ ${\cal U}$ is an open cover of $M$, then there exists an
aperiodic oriented 1-foliation of $M$ in the same smoothness category, such
that each leaf is contained in an element of \  ${\cal U}$, and
whose minimal sets have codimension one.  \label{thgreg} \end{theorem}

\begin{theorem} Let $M$ be a $\PL$ manifold of dimension $n\geq 3$,  $1\leq
k\leq n-1$, and let \ ${\cal U}$ be an open cover of $M$.  A $\PL$  1-foliation
of $M$ can be modified in a $\PL$ fashion so that each leaf is contained
in an element of \ ${\cal U}$, there are no circular leaves, and all minimal sets are
$k$-dimensional.  \label{thgregs} \end{theorem}

The stopped set of each of the plugs constructed in this paper contains an arc. If
the manifold $M$ is closed, then the foliation specified in Wilson's theorem has
finitely many minimal sets. A plug can be ``weaved'' through the foliation so the
arc in the stopped set meets a dense leaf in each of the minimal sets. Hence for closed
manifolds, we get the following variations of the above theorems:

\begin{theorem} If $M$ is a continuous, $C^r$, $C^\infty$, $C^\omega$, or $\PL$ closed manifold
of dimension $\geq 3$ admitting an oriented 1-foliation in the same smoothness
category,  then there exists an
aperiodic oriented 1-foliation of $M$ in the same smoothness category, with exactly one
minimal set, and the minimal set has codimension one. \label{thcompact}
\end{theorem}

\begin{theorem} Let $M$ be a closed $\PL$ manifold of dimension $n\geq 3$ and $1\leq
k\leq n-1$.  A $\PL$  1-foliation of $M$ can be modified in a $\PL$ fashion so that
there are no circular leaves, and  there is exactly one minimal set which is 
$k$-dimensional.  \label{thcompactPL} \end{theorem}

\section{\label{sinsert} Self-insertion}

As a warm-up to the important technique of self-insertion, we first
define and analyze partial insertion of plugs.

\fpartial

Let ${\cal P}$ be an $n$-dimensional plug with base $F$ and support $P$. Let
$\alpha:N_{F \times I} \to N_P$ be an attaching map for ${\cal P}$. Let $X$ be an
$n$-manifold with an oriented $1$-foliation ${\cal X}$.   Let $D$ be a closed
domain in $F$ whose boundary in $F$ (in the point-set topology sense of
boundary) is a properly embedded $(n-2)$-manifold $A$ transverse to $\d F$.  It
may happen that either $A$ or $D$ is not connected. The domain $D$ is, in
particular, a manifold with corners.  A {\em partial insertion map\/} of $\cal
P$ is an embedding $\sigma:D \to X$ such that $\sigma(A) \subset \d X$, but
$\sigma(D - A)$ is disjoint from $\d X$.  The map $\sigma$ can be extended to an
embedding of $D \times I$ which takes vertical fibers to leaves. Define
$\hat{X}$ to be $X-\sigma((D\times I)-N_{F \times I})$ glued to $P$ by the
composition $\sigma \circ \alpha^{-1}$ where it is defined.  Define $\hat{\cal
X}$ by gluing ${\cal X}$ to $\cal P$ in the same way.

Figure~\ref{fpartial} shows an example of a partial insertion of a plug with
base $D^2$.  The partial insertion creates transverse boundary (resembling a
ledge and an overhang); unlike with a complete insertion, some leaves in the
result have endpoints.  By the matched ends condition, all leaves of the partial
insertion either have both endpoints or neither endpoint at the transverse
boundary.  A leaf $\hat{{\cal X}}$ from the depths of $X$ might enter $P$ and
never return to $X$, but only by approaching some minimal set of ${\cal P}$ and
never by terminating.

Let ${\cal P}$, $A$, and $D$ be as before.  A {\em self-insertion map\/} for $\cal
P$ is an embedding $\sigma:D \to P$ such that the image of $\sigma$ is
transverse to ${\cal P}$ and disjoint from $F_\pm$,  $\sigma(A) \subset
\d P$,  $\sigma(D-A)$ is disjoint from $\d P$,  $\sigma(A)$
avoids leaves with endpoints in $\alpha _-(A)$, and $\sigma(A)$ does not intersect any leaf
twice.  Choose an attaching map $\alpha:N_{F \times I} \to N_P$ such that the image of
$\sigma$ is disjoint from the closure of $\alpha((D \times I) \cap N_{F \times
I})$.  The map $\sigma$ can be extended to an embedding of $D \times I$ which is
also disjoint from the closure of $\alpha((D \times I) \cap N_{F \times I})$ and
from $F_\pm$.  Let $\hat{P}$ be the space $P - \sigma((D \times I)- N_{F \times I})$
identified to itself by $\sigma \circ \alpha^{-1}$ where it is defined. The
self-insertion of ${\cal P}$ at $\sigma$ is the foliation $\hat{\cal P}$ on
$\hat{P}$ obtained by the same gluing.

Figure~\ref{fentryexit} shows an example of self-insertion of an un-plug with
base $I$.  Keeping the above notation, the surfaces $\alpha_-(D-\d (A))$ and
$\alpha_+(D-\d (A))$, which are on the boundary of $P$, are in the interior of
$\hat{P}$.  These are the {\em internal entry\/} and {\em exit regions\/}, respectively, of
$\hat{{\cal P}}$.  Similarly, $F_- - \alpha_-(D-\d (A))$ and $F_+ - \alpha_+(D-\d
(A))$ are the {\em external entry\/} and {\em exit regions\/} of $\hat{{\cal P}}$. Many leaves
of $\hat{{\cal P}}$ cross the internal entry and exit regions once or many times. 
The points where a leaf $l$ does so are the {\em transition\/} points of $l$, or the
entry and exit points. In addition, $l$ may begin with an external entry and/or
end with an external exit.  Progressing along a leaf $l$ of $\hat{{\cal P}}$ in
the positive direction, we encounter a history of transition points. In
general $l$ consists of segments of leaves of ${\cal P}$ separated by
transition points. If $p \in l$ is an entry (exit), the leaf of ${\cal P}$
preceding $p$ is {\em interrupted\/} ({\em interrupting\/}) leaf $p$, while
the leaf following $p$ is the interrupting (interrupted) leaf.

In the self-insertion in Figure~\ref{fentryexit}a,
a leaf has the following history of transition points:
\begin{center}
(external) entry, entry, entry, exit, exit, (external) exit.
\end{center}
In the self-insertion in Figure~\ref{fentryexit}b, a leaf entering the bottom
has the history
\begin{center} (external) entry, entry, entry, entry, $\ldots$ \end{center}
In a given history, an entry and an exit are {\em matched\/} if the exit
follows the entry and they satisfy the following inductive rule:  They are
matched if they are adjacent, and otherwise they are matched if all entries
and exits between them are matched to each other.  For example, if a leaf in
a hypothetical self-insertion has the history given in Figure~\ref{frainbow},
then the entries and exit are matched as indicated.  The proof of the
following lemma is an instructive exercise for the reader; it can also be
found in \cite{kseifert}.  Note that the lemma would not hold for a
hypothetical definition of self-insertion of a semi-plug.

\frainbow

\begin{lemma} If two transition points are matched, they  have the same
interrupted leaf and they are the endpoints of the same interrupting leaf.
\label{linterrupt}
\end{lemma}

The definition of matched transition points and Lemma~\ref{linterrupt} suggest
an interpretation of the geometry of a self-insertion in terms of a recursive
algorithm.  The following procedure follows a leaf of $\hat{{\cal P}}$ by
following segments of leaves of ${\cal P}$.  The input to the procedure is a base
point $p \in F$; the procedure follows the corresponding leaf $\hat{l}$ of
$\hat{{\cal P}}$.

\fentryexit

\begin{description}
\item[\,]Procedure {\sc followleaf}($p$)
\item[1.]Let $l$ be the leaf of ${\cal P}$ that begins at $\alpha_-(p)$.
\item[2.]Follow $l$ until some $\sigma(q)$ or $\alpha_+(p)$ is reached.
\item[3.]If $\alpha_+(p)$ is reached, quit.
\item[4.]If $\sigma(q)$ is reached, then do {\sc followleaf}($q$) and go to step 2.
\end{description}

In computer programming terms, a recursive procedure such as {\sc followleaf} is
implemented by means of a pushdown {\em stack\/}:  When {\sc followleaf} is
called, its argument is pushed onto the top of the stack.  It is removed when
that particular instantiation of {\sc followleaf} terminates.  Given a point $p
\in \hat{{\cal P}}$, the stack of $q$ is the sequence of points $E(p)$,
$E(\sigma(E(p))$, $E(\sigma(E(\sigma(E(p)))),\ldots,$ where $q = E(p)$ is a
point in $F$ such that at least one of $\alpha_+(q)$ and $\alpha_-(q)$ is an
endpoint of the leaf of ${\cal P}$ containing $p$.  The function $E$ is not
defined on all of $P$; the stack at $q$ extends as long as $E$ is defined and
may therefore be either a finite or infinite sequence.  It is easy to check that
as $q$ moves along its leaf, its stack changes in the same way as the stack used
to execute {\sc followleaf}.

Inspection of {\sc followleaf} demonstrates that in a finite leaf, all
transition points are matched.  In particular, the two endpoints are matched to
each other.

This can be construed as saying that $\hat{{\cal P}}$ has matched ends. 
Unfortunately, $\hat{{\cal P}}$ is not even a flow bordism, because the top and
bottom of $\hat{P}$ have stair-steps and are not transverse to $\hat{{\cal P}}$. 
However, without disturbing the transition point scheme, it is possible to
approximate the top and bottom of $\hat{P}$ by surfaces $\tilde{F}_+$ and
$\tilde{F}_-$ that are transverse to $\hat{{\cal P}}$.  Let $\tilde{P}$ be the
manifold that results from cutting along these surfaces, and let $\tilde{\cal
P}$, the self-inserted plug, be the restriction of $\hat{{\cal P}}$ to $\tilde{P}$.

As Figure~\ref{fentryexit}b demonstrates, a leaf can have unmatched transition
points.  Indeed, by Lemma~\ref{linterrupt}, if a transition point (whether
internal or external) is the endpoint of an infinite leaf of ${\cal P}$, then it
cannot be matched in $\hat{{\cal P}}$.  Its leaf in $\hat{\cal P}$ (and hence in
$\tilde{{\cal P}}$) is therefore either infinite or circular, and by
Lemma~\ref{lisplug}, $\tilde{{\cal P}}$ is a plug. In particular, if $p \in F$ is
in the stopped set of ${\cal P}$ and $\alpha _-(p)$, then it lies in an infinite
leaf of $\hat{{\cal P}}$. Thus, a self-insertion cannot close infinite leaves of
this type.  On the other hand, Figure~\ref{fentryexit} also demonstrates that a
self-insertion can create new circular leaves; these necessarily have unmatched
transition points.

Finally, since a self-insertion changes the base of a plug, it may happen that
the self-insertion of an insertible plug is not insertible. It may also happen
that the self-insertion of an untwisted plug is twisted.

The following example shows that a self-insertion may cause the boundary of a
plug to disappear completely. Let ${{\cal Q}}$ be a plug with support $D^2\times I$
such that  the transverse boundary consists of the top and bottom of the
cylinder, and the side of the cylinder is the parallel boundary. A
self-insertion of ${{\cal Q}}$ using an immersion of the entire disk $D^2$ 
results in a plug with a boundaryless support ({\em e.g.\/}, $S^2\times S^1$). 

\section{\label{scounter} An analytic plug}

This section contains a construction, similar to that of \cite{kseifert}, of a
two-component self-insertion of a plug ${\cal W}$ which is a concatenation of two
semi-plugs ${\cal W}_s$ and $\bar{\cal W}_s$, each a mirror image of the other.

\fwilson

Parametrize $S^1$ by $\theta$ with $0\leq \theta <10$.  Let  be $F=[-1,1]\times
S^1$ ($F$ is an annulus), and the support of ${\cal W}_s$ be the cylinder $W_s = F
\times [-1,1]$.  The coordinates of $W_s$ are $r$, $\theta$, and $z$, similar to
the  cylindrical coordinates in ${\Bbb R}^3$. The support of $\bar{{\cal W}}_s$ is
also $F \times [-1,1]$, but it will be denoted $\bar{W}_s$ to avoid confusion;
let $W = W_s \cup \bar{W}_s$.

The semi-plugs ${\cal W}_s$ and $\bar{\cal W}_s$ are generated by the vector fields
$$\vec{W_s}(r,\theta,z) = {\d\over\d\theta} + (r^2+z^6){\d\over\d z} $$
and 
$$\bar{\vec{W}}_s(r,\theta,z) = -{\d\over\d\theta} + (r^2+z^6){\d\over\d z}\, .$$
The concatenation identifies the top $F \times \{1\}$ of $W_s$ with the bottom
$F \times \{-1\}$ of $\bar{W}_s$; the map $\alpha_-:F\to W$ sends $F$ onto the
bottom of $W$ by $\alpha_-(p) =(p,-1)\in W_s$.  The vector fields  parallel to
${\cal W}_s$ and $\bar{\cal W}_s$  are oriented in the positive $z$ direction.
Since $W_s$ and $\bar{W}_s$ are trivially foliated in the respective
neighborhoods of $F\times \{1\}$ and $F\times \{-1\}$ and they
are analytic, the plug ${{\cal W}}$ is analytic.

Let $T=\{(r,\theta ,z)\in W_s \st  r=0, z=0\}$ and  $\bar{T}=\{(r,\theta
,z)\in \bar{W}_s\st  r=0, z=0\}$ be the circles of ${\cal W}_s$ and
$\bar{{\cal W}}_s$.

We define a self-insertion whose  domain $D\subset F$ consists of two components
$D_s$ and $\bar{D}_s$. The self-insertion map $\sigma:D \to W$ divides into a map
$\sigma:D_s \to W_s$ given by the formula
$$\sigma(r,\theta) = (r-\frac14r^2-2(\theta-2)^2,6\, ,2-\theta ),$$
where $D_s$ is the set of all points $(r,\theta)$ such that 
$(r-\frac14r^2-2(\theta-2)^2,6\, ,2-\theta)$ lies in $W_s$, and a map
$\sigma:\bar{D}_s \to \bar{W}_s$ similarly defined by
$$
\sigma(r,\theta) = (r-\frac14r^2-2(\theta-8)^2,4\, ,\theta -8),
$$
where $\bar{D}_s$ is the set of all points $(r,\theta)$ such that 
$(r-\frac14r^2-2(\theta-8)^2,4\, ,\theta -8)$ lies in $\bar{W}_s$. 

To check that the equations for $\sigma$ form a valid self-insertion, we first
note that $D_s$ and $\bar{D}_s$ are disjoint, because if $(r,\theta)\in D_s$,
then  $1 < \theta < 3$, while if $(r,\theta)\in \bar{D}_s$, then $7 <
\theta < 9$.  Each leaf of $\d W$ intersects $\im(\sigma)$ at most once,
because the slope ${dz \over d\theta}$ of any leaf in $\d W$ is bigger than $1$
and is positive in $W_s$. The variable $\theta$ would have to increase by 8 for
a leaf to  connected $\sigma(D_s)$ with $\sigma(\bar{D}_s)$ and by 10 for it to
connect either component with itself, but this is impossible. Finally, the
condition that $\sigma$ avoid leaves in $\d W$ that begin in its domain is
automatically satisfied, because such leaves have $r = 1$ while $\im(\sigma)$
only intersects $\d W$ at $r = -1$. The self-inserted plug $\tilde{{\cal W}}$
given by $\sigma$ is pictured in 
Figure~\ref{fdad}.

\fdad

Observe that $\sigma(0,2) \in T$ and that $\sigma(0,8) \in \bar{T}$ and that
$(0,2)$ and $(0,8)$ are both in the stopped set of ${\cal W}$. In other words, the
self-insertion breaks both circles of ${\cal W}$.  Moreover, if $(r,\theta,z) =
\sigma(r',\theta')$, then $r \le r'$, with equality occurring only for two
points of $F$; one point is sent to the circle $T$ by $\sigma$ and the other is
sent to $\bar{T}$.  This is the important {\em radius inequality\/} for
self-insertions.

A point $p \in \tilde{W}$ may be considered as a point in $W$.  To make a
distinction between the leaves of $\tilde{{\cal W}}$ and $\cal W$ containing $p$, the leaves
are denoted by $\tilde{l}$ and $l$.

\begin{theorem} The self-inserted plug $\tilde{{\cal W}}$ has no circular
leaves. \label{thnocirc}
\end{theorem}
\begin{proof} Let $p_1,p_2,\ldots$ be the stack of some point in $\tilde{W}$. 
By the radius inequality, the points have strictly decreasing radii $r_1 >
r_2 > \ldots \,.$   If some $p_n$ is either $(0,2)$ or $(0,8)$, then the
stack ends at $p_n$, because $\sigma(p_n)$ lies on the circle $T$ or
$\bar{T}$ and $E(\sigma(p_n))$ is undefined. At all other points of $D_s\cup
\bar{D}_s$, the radius inequality is a strict inequality.

Suppose that $\tilde{l}$ is a circular leaf of $\tilde{{\cal W}}$ and let $p$ be
a point that varies along $\tilde{l}$. As $p$ goes all the way around
$\tilde{l}$, the stack of $p$  either grows, shrinks, or stays at the same
level.  If it grows or shrinks, then it must be infinite, and moreover it
must be eventually periodic, which is impossible by the radius inequality. 
If it stays at the same level, then for some $q \in l$ the stack reaches a
minimum height; the segment of $\tilde{l}$ containing $q$ must therefore
correspond to some leaf $l$ of ${\cal W}$ which is interrupted but does not
interrupt other leaves.  The leaf $l$ must be a circle, and must therefore be
one of the two circles $\tilde{T}$ and $\tilde{\bar{T}}$.  However, the
transition points on these circles lead to the entry points $(0,2)$ and
$(0,8)$ which lie in the stopped set of ${\cal W}$ and cannot be matched.  In
conclusion, all avenues for a circular leaf in $\tilde{{\cal W}}$ lead to
contradiction.
\end{proof}

\begin{lemma} If $\tilde{l}$ is a leaf of $\tilde{{\cal W}}$ whose radii avoid the
interval $(-\epsilon,+\epsilon)$ for some $\epsilon >0$, then $\tilde{l}$ is a
finite leaf.  \label{lrlimit}
\end{lemma}
\begin{proof}
Let $r_1 > r_2 > \ldots$ be the radii of the stack of some point on such a
leaf $\tilde{l}$. By hypothesis the points of the stack are at least
$\epsilon$ away from the critical points $(0,2)$ and $(0,8)$.  Therefore
there exists a $C$ such that $r_n > C + r_{n+1}$, which in turn implies that
the height of the stack is bounded by $2/C$. Moreover, all leaves of ${\cal W}$
that intersect the image of $\sigma$ infinitely many times have $r=0$, so
there exists an $N$ such that all component leaves of $\tilde{l}$ intersect
$\sigma$ at most $N$ times.  Modelling the behavior of $\tilde{l}$ by {\sc
followleaf}, each instantiation of {\sc followleaf} can only call {\sc
followleaf} at most N times, and the recursion can only extend to a depth of
$2/C$.  Therefore there are at most $N^{2/C}$ internal calls, which means
that $\tilde{l}$ has at most $N^{2/C}$ internal entry points, or $2N^{2/C}$
transition points in total.  This implies that $\tilde{l}$ is finite, since
none of the component leaves is infinite. 
\end{proof}

\begin{lemma} If $\tilde{l}$ is a leaf of $\tilde{{\cal W}}$ with entry point at
$r=0$, then all subsequent transition points are matched and the positive limit
set of $\tilde{l}$ contains the (formerly circular) leaf $\tilde{T}$.
\label{lrbig}
\end{lemma}
\begin{proof} At the bottom level, the leaf $\tilde{l}$ follows a leaf $l$ of
${\cal W}$ which approaches $T$.  This leaf $\tilde{l}$ meets an infinite sequence
of transition points, all of which have $r > 0$. Since the radii of higher-level
transition points are only larger, Lemma~\ref{lrlimit} implies that all of these
transition points are matched.
\end{proof}

A similar argument shows that the negative limit set of a leaf of $\tilde{\cal
W}$ with exit point at $r=0$ contains the leaf $\tilde{\bar{T}}$.

\begin{theorem} The plug $\tilde{{\cal W}}$ has a unique non-trivial minimal set
and it is 2-dimensional. \label{thdim}
\end{theorem}
\begin{proof}
By Lemmas~\ref{lrlimit} and \ref{lrbig}, the closure of any infinite leaf
contains the leaf $\tilde{T}$.  Therefore there is only one minimal set and it
contains $\tilde{T}$. Since $\sigma(0,2) \in T$, to understand the minimal set
it suffices to follow the leaf $\tilde{l}$ of $\tilde{{\cal W}}$ starting at
$(0,2,-1) = \alpha_-(0,2)$.

Let $A(x)$ be the antiderivative of $1/(1+x^6)$ with $A(0) = 0$; in particular,
$$\int_{-\infty}^{\infty} {1 \over 1+x^6} dx = 2A(\infty)$$
by abuse of notation.  Solving the differential equation defined by $\vec{W}$,
the leaves of ${\cal W}_s$ are given by the equation
$$\theta = -{1 \over 5z^5} + C$$
when $r = 0$, and by
$$\theta = r^{-5/3} A(r^{-1/3}z) + C$$
when $r \ne 0$.  If $l$ is the leaf of ${\cal W}_s$ starting at $(0,2,-1)$, then
$l$ is given by the former equation with $C = 9/5$.  For each integer $n \ge
0$, $l$ intersects the surface $\theta = 6$ with
$$z = z_n = - (21 + 50n)^{-1/5},$$
and for $n$ large enough, these points also lie in $\sigma(D_s)$. For such $n$,
choose $(r_n,\theta_n) \in F$ so that
$$\sigma(r_n,\theta_n) = (0,6,z_n).$$
Then the numbers $\theta_n = 2-z_n  \to 2$ as $n \to \infty$, while $r_n
\approx 2 z_n^2$ as $n \to \infty$ by inspection of the relation
$$r_n - \frac14 r_n^2 - 2 z_n^2 = 0.$$
Let $l_n$ be the leaf of ${\cal W}_s$ beginning at $(r_n,\theta_n,-1)$.  By
Lemma~\ref{lrbig}, in the leaf $\tilde{l}$, each leaf $l_n$ concatenated with
its mirror image interrupts the leaf $l$. Let $(r_n,\theta'_n,1)$ be the other
endpoint of\ $l_n$ in $W_s$; $\theta'_n$ is given by
$$\theta'_n = \theta_n + 2r_n^{-5/3} A(r_n^{-1/3}).$$
Combining several identities and approximations,
$$\theta'_n \approx 2+2^{1/3}(21+50n)^{2/3}A(\infty)$$
for $n$ large.  Thus, $\theta'_n \to \infty$, but  $\theta'_{n+1} - \theta'_n \to
0$ as $n \to \infty$.  It follows that $\{\theta'_n \bmod 10\}$ is a dense
subset of the circle $\R/10\Z$, and that the circle $K$ of all points
$(0,\theta,1) \in W_s$ is in the closure of $\tilde{l}$.  Following the
leaves of ${{\cal W}}_s$ and $\bar{\cal W}_s$, the circle $K$ sweeps out the
entire annulus $N$ between $T$ and $\bar{T}$. The surface $N$, except where
it is cut by $\im(\sigma)$, is therefore in the closure of $\tilde{l}$ in
$\tilde{{\cal W}}$, which demonstrates that the minimal set of $\tilde{\cal W}$
is 2-dimensional.
\end{proof}

Finally, the following proposition was  pointed out to the authors by
\'E.~Ghys, who credits S.~Matsumoto \cite{Ghys}.  It slightly simplifies
Section~\ref{sstoppage}, since it means that it is not necessary to insert
$\tilde{{\cal W}}$ into another content-stopping plug.

\begin{proposition} The plug $\tilde{{\cal W}}$ stops content. \label{pcontent}
\end{proposition}
\begin{proof} The set of points $(r,6,z) = \sigma(r',\theta) \subset W_s$ with 
$r' \le 0$ is the region in the $\theta = 6$ plane with
$r \le -2z^2$.  In this region, the slope of a leaf of ${\cal W}$ satisfies
$${dz \over \d \theta} = r^2 + z^6 \le \sqrt{2}r^2.$$
If $l$ is a leaf of ${{\cal W}}$ at radius $r$ with $0 > r > -\frac 1{20}$, then
in the region $r \le -2z^2$, the successive intersections of $l$ with
$\im(\sigma)$ are less than $10\sqrt{2}r^2 < \frac r{10}$ apart.  But the
region $r \le -2z^2$ has width $\sqrt{2r} > \frac r{10}$ in the
$z$ direction, so it follows that
the leaf in $\tilde{l}$ in $\tilde{{\cal W}}$ containing segments of $l$
eventually meets an internal entry point with interrupting radius $r'$ such
that $0 \ge r' > r$. All entry points that $l$ meets at radius $r' > 0$ are
matched by Lemma~\ref{lrbig}, but an entry at radius $0$ cannot be matched,
and an entry at radius $0 > r' > r$ begins a new leaf $l'$ of ${\cal W}$ with
the same properties as $l$.  Therefore the stack of $\tilde{l}$ grows
indefinitely and $\tilde{l}$ is an infinite leaf.  In conclusion, all points
of $\hat{F}$ with radius $0 > r > -\frac 1{20}$ are in the stopped set of
$\tilde{{\cal W}}$.
\end{proof}

Proposition~\ref{pcontent} applies equally well to any conceivable variation
of the formulas defining the plug $\tilde{{\cal W}}$, with the conclusion 
that no such analytic (or even $C^1$) plug can preserve volume given by a
volume form.  On the other hand, Proposition~\ref{pcontent} is valid for some
but not all PL self-insertions.  Nevertheless, it seems impossible to satisfy
the radius inequality and preserve volume also.  Reference \cite{kvol} gives
a version of the Schweitzer plug which preserves volume.  Like
Schweitzer's example, it is $C^1$.

\section{\label{shigher}Higher dimensions}

We first consider oriented 1-foliations in $n$ dimensions with $n > 3$.

Let $T^{n-2}$ be an $(n-2)$-dimensional torus, parametrized by coordinates
$\theta_1, \theta_2, \ldots, \theta_{n-2}$ with period $10$.  Let
$$
\vec{\theta} = \sum_i k_i\frac{\d}{\d\theta_i}
$$ 
be a vector field, where the coefficients $\{k_i\}$ are between $0$ and $1$ and
are linearly independent over $\Q$, the rational numbers.  All leaves of the
foliation parallel to $\vec{\theta}$ are dense.  Parametrize the manifold
$W_{s,n} = [-1,1] \times T^{n-2} \times [-1,1]$ by the coordinates
$r,\theta_1,\ldots,\theta_{n-2},z$, and define the vector field $\vec{W}_{s,n}$
on $W_{s,n}$ by the formula
$$\vec{W}_{s,n} = \vec{\theta} + (r^2+z^6){\d\over\d z}.$$
The foliation ${\cal W}_{s,n}$ parallel to $\vec{W}_{s,n}$ is a semi-plug, and the
mirror-image construction yields a plug ${\cal W}_n$ with base $I \times T^{n-2}$.
The plug ${\cal W}_n$ is aperiodic and is similar to a construction of Wilson
which settled the higher-dimensional Seifert conjecture for all manifolds with
Euler characteristic 0.  However, its minimal sets are at most $(n-2)$-dimensional; to
achieve $(n-1)$-dimensional minimal sets we will perform a self-insertion.

Let $a = (2,2,,\ldots,2)$ in the parametrization by the $\theta_i$'s, and let
$b = (8,8,\ldots,8)$. Consider a self-insertion with two components, one
component going into ${\cal W}_{s,n}$ and the other going into $\bar{\cal
W}_{s,n}$. The formula for the part of $\sigma$ that maps into $W_{s,n}$ is
$$\sigma(r,\theta_1,\ldots ,\theta_{n-2}) =
(r-\frac14r^2-2\sum _1^{n-2}(\theta_i-2)^2,6\,\ldots ,6\, 
,2-\theta_1+\sum _2^{n-2}(\theta_i-2)),$$
while the formula that maps into $\bar{W}_{s,n}$ is
$$\sigma(r,\theta_1,\ldots ,\theta_{n-2}) =
(r-\frac14r^2-2\sum _1^{n-2}(\theta_i-8)^2,4\,\ldots ,4\, 
,\sum _1^{n-2}(\theta_i-8)).$$
Let $\tilde{{\cal W}}_n$ be the self-inserted plug.

The self-insertion clearly satisfies the radius inequality with respect to the
variable $r$.  The base after self-insertion is $I \times T^{n-2}$ with the
boundary components connected by two boundary handles, which admits a  bridge
immersion in $\R^{n-1}$ as illustrated in Figure~\ref{ftoungue} in the case
$n=4$.  Following the argument for three dimensions, if $l$ is the leaf of
${\cal W}_s$ beginning at $\alpha_-(0,2,2,\ldots,2)$, then $\tilde{l}$ is in
the minimal set of $\tilde{{\cal W}}_n$.  In $\tilde{l}$, the leaf $l$ is
interrupted by a sequence of leaves $l_n$ with radius $r_n \to 0$.  In the
$\theta_i$ directions, these leaves progress in the $\theta_i$ directions by
an amount that goes to infinity as $n \to \infty$, but at a rate that goes to
0.  It follows that the entire surface with $r = 0$ and $z=1$ in $W_{s,n}$ is
in the closure of the leaves $l_n$, which in turn forces the minimal set of
$\tilde{{\cal W}}_n$ to have codimension 1.

\ftoungue

In the general case, the goal is to open the leaves of a $k$-dimensional foliation
of an $n$-manifold so that the closure of any such leaf is at least
$(n-1)$-dimensional.  If $k = n-1$, there is almost nothing to prove.  To
eliminate the possibility of a minimal set of codimension 0 (the whole
manifold), choose a circle or properly embedded line transverse to the foliation 
and put in a Reeb structure along this curve.

The interesting case is $1 < k < n-1$.  Following Schweitzer
\cite{Schweitzer}, consider the manifold $W_{n-k+1} \times S^k$ foliated by
leaves $l \times S^k$, where $l$ is a leaf of ${\cal W}_{n-k+1}$ and $S^k$ is
the $k$-sphere. The techniques of insertion and  global stoppage readily
generalize to such a foliation, and it is easy to check that the only 
minimal set has codimension 1.

\section{\label{spl}The PL case}

The PL construction mostly uses the same geometric ideas as the analytic (and
therefore smooth) case, but they are realized somewhat differently because PL
foliations do not have well-behaved parallel vector fields. The construction
has one new detail, a vertical annulus of circles, which gives the minimal
set high dimension in a different way.

Let $H$ be a compact manifold, usually with boundary or corners.  Let $f:H
\times [a,b] \to H \times [a,b]$ be a PL homeomorphism and let $0<l<1$ be a real
number.  Let ${\cal L}$ be the foliation of $H \times [a,b] \times [0,1]$ such that,
for fixed $p \in H$ and $z \in \R$, the set $\{(p,z + lx,x) \st  x\in [0,1],z+lx \in
[a,b]\}$ is a leaf.  Orient all such leaves from $H \times [a,b] \times \{0\}$ to
$H \times [a,b] \times \{1\}$. The {\em slanted suspension} of $f$ with {\em slant}
$l$ is defined as $Z$, the manifold $H \times [a,b]\times [0,1]$ with
$(f(p,z),0)$ identified with $(p,z,1)$, together with the foliation ${\cal Z}$
induced from ${\cal L}$.  The slanted suspension $\cal Z$ is a PL foliation by
construction, and moreover is a flow bordism with entry and exit regions $H
\times S^1$.

\fcollar

Let $H = [-1,1]$ and $[a,b]=[-2,2]$.  Let $f:[-1,1]\times [-2,2] \to [-1,1]\times
[-2,2]$ be a PL homeomorphism which is the identity at $\d([-1,1]\times [-2,2])$,
such that $f(0,0) = (0,-1)$, $f(0,1) = (0,0)$, and such that $f$ is linear on the
two triangles and two trapezoids illustrated in  Figure~\ref{fcollar}.  It is easy
to check that the slanted suspension of $f$ with slant $l=1$ is a semi-plug with an
annulus of circular leaves corresponding to the line segment from $(0,-1)$ to
$(0,0)$. It is  convenient to take a 20-fold covering of the slanted
suspension, to obtain a semi-plug ${{\cal W}}_{PL,s}$ whose support 
$W_{PL,s}=[-1,1]\times [0,20]\times [-2,2]$ (where 0=20) is parametrized 
by $r$, $\theta$,  and $z$ with the suspension direction $\theta$.
 The semi-plug ${\cal W}_{PL,s}$ is
similar to the analytic semi-plug ${\cal W}_s$, except that it has an annulus of
circular leaves instead of one circular leaf. Let $\bar{{\cal W}}_{PL,s}$ (with support
$\bar{W}_{PL,s}$) be a the  mirror image of ${\cal W}_{PL,s}$ obtained by changing $z$ to
$-z$. Let ${\cal W}_{PL}$ (with support
${W}_{PL}$) be the plug obtained from ${\cal W}_{PL,s}$ by the mirror-image
construction, the concatenation of ${\cal W}_{PL,s}$ and $\bar{\cal W}_{PL,s}$. 

Although ${\cal W}_{PL}$ has two annuli of circles, $T$ and $\bar{T}$, it is still possible
to break all of them with a self-insertion, since the stopped set is also a circle.  The
radius inequality preserving self-insertion $\sigma :D\to W$  is defined on two parts of
$D$:  $D_s$ containing the segment with endpoints $(0,9)$, $(0,10)$ which is mapped
linearly onto the whole segment of  $T$ at $\theta =20$, and $\bar{D}_s$ contaning the
segment with endpoints $(0,14)$, $(0,15)$ which is mapped linearly onto the whole segment
of  $\bar{T}$ at $\theta =5$. (Figure~\ref{fvertical} shows the $D_s$ 
part of $\sigma$.)

\fvertical

Its geometry is more or less the same as that in the
analytic case, since the leaves of ${\cal W}_{PL}$ are  at constant $r$ and the
self-insertion satisfies the radius inequality. Lemma~\ref{lrlimit}, slightly modified
Lemma~\ref{lrbig}, and Theorem~\ref{thdim}  hold for $\tilde{{\cal W}}_{PL}$. The
internal entries  of a point at $r = 0$ have radii that converge to $r = 0$, and their
leaves converge to the annuli of circles of ${\cal W}_{PL}$. Therefore the remnants of
these annuli after self-insertion are contained in the minimal set, and the minimal set
is still 2-dimensional.

The construction for higher-dimensional foliations also apply to the PL case
without modification.

\section{\label{ssymbolic}Symbolic Dynamics}

Although most aperiodic self-inserted Wilson-type plugs have 2-dimensional
minimal sets, a carefully chosen self-insertion may result in a 1-dimensional
minimal set.  Such a self-insertion has interesting symbolic dynamics which can
be described explicitly. The self-insertion is easiest to define in the 
continuous category, but with yet more care it is also possible in the PL
category.

\begin{theorem} There exists an aperiodic $\PL$ plug with 1-dimensional
minimal sets. \label{th1dplug}
\end{theorem}

\begin{proof}
First, we consider a self-insertion that only breaks one circle.  The
other circle is the unique minimal set, but the leaf containing the broken
circle is in a 1-dimensional invariant set which resembles the minimal set in
the final construction.

Let $R = [-1,\frac52] \times [-3,\frac32]$ be a rectangle parametrized by $r$ and
$z$, and let $f:R \to R$ be the PL homeomorphism illustrated in
Figure~\ref{fmapf}.  The map $f$ is linear in each of the regions delineated by
the four rays $z = -r \le 0$, $z = -r \ge 0$, $z = \frac r2 \ge 0$, and $z = 2r
\le 0$.  In the upper region, $f(r,z) = (r,2z-\frac32)$; in the lower region,
$f(r,z) = (r,\frac z2-\frac32)$,  and in the two side regions, $f$ matches the
unique linear transformation which makes it continuous along the four rays. Note
that $f(0,0)=(0,-\frac32)$, and for any other point $(r,z)\in R$, the difference
between $z$ and the $z$-coordinate of $f(r,z)$ is greater than $-\frac32$. Let
${\cal V}_s$ be the slanted suspension of $f$ with slant $\frac32$, and let $V_s$ be
its support. More precisely, $V_s$ is obtained from $[-1,\frac52]\times [0,1]
\times [-3,\frac32]$ by identifying  $(\bar{r},0,\bar{z})$ with $(r,1,z)$, where
$(\bar{r},\bar{z})=f(r,z)$. We use the coordinates $(r,\theta ,z)$, and the  suspension
is in the direction of $\theta$. Each leaf of ${\cal V}_s$ is the  union of segments
$\{(r,\theta ,z + {\frac32}\theta) \st  r\in [-1,\frac52], z \in \R, \theta \in [0,1],
z+{\frac32}\theta  \in [-2,2]\}$. The semi-plug ${\cal V}_s$ is another PL
analogue of the analytic semi-plug ${\cal W}_s$. Like $\cal W_s$, it has one circular
leaf $T$. Note that $T$ passes through the point $(0,0, -\frac32 )=(0,1,0)$.  

\fmapf

Let $\bar{{\cal V}}_s$ be the mirror image of ${\cal V_s}$, and let $\bar T$ be the
circular leaf of $\bar{{\cal V}}_s$. Let ${\cal V}$ be the concatenation of the two
mirror-image semi-plugs ${{\cal V}_s}$ and $\bar{\cal V}_s$. The support $V$ of ${\cal
V}$ is obtained by idetifying the top of ${ V_s}$ with the bottom of $\bar{ V}_s$. Let
$F$ be the  base of ${\cal V}$.  

Let $B=\{(r,\theta ) \st  r\geq 0, \frac r{12} \le {\theta - \frac13} \le
\frac r4\}$.  We define a self-insertion map $\sigma:D \to V_s \subset V$ on an
appropriate   $D\subset   F$ containing $B$.  The image of $\sigma$ lies
in the $\theta = 1$ section of $V_s$, and 
 $\sigma $  is described    explicitly on $B$ by:  
 $$g(r,\theta) = (\frac12r,9\theta - \frac74 r - 3)=(x,y),$$ 

 $$h(x,y) = \left\{\begin{array}{lr} 
(2x-2y,1,\frac x2) & \mbox{if $ y \ge \frac12 x$}\\ 
(x,1,y) & \mbox{if $-x < y < \frac12 x$}  \\
(2x+y,1,-x) & \mbox{if $ y \le -x $} 
\end{array}\right.,$$
and  $\sigma(r,\theta) = h \circ g$.   On $D-B$,  $\sigma$  tapers in the $z$
direction to have $\im(\sigma) \cap \d V_s$  small enough so that leaves of $\d V_s$
intersect $\im(\sigma)$ at most once, and none of the leaves that do intersect have an
endpoint in $D$.  Figure~\ref{fsigma} shows the shape of
$\im(\sigma)$ in the section of $V_s$ with $\theta = 1$.  The formulas
demonstrate that $\sigma$ satisfies the radius inequality on $B$. Outside $B$ the
radius inequality is easy to achieve.  Let $\tilde{{\cal V}}$ be $\cal V$
self-inserted by $\sigma$. Let $\tilde T$ be the leaf of $\tilde{{\cal V}}$ containing
segments of $T$. By Lemma~\ref{lrbig}, $\sigma |_B$ determines the geometry of
the leaf $\tilde T$.

\fsigma

Let $E_1,E_2,\ldots$ be a sequence of closed disks in $B$ such that
$$E_n=\{(r,\theta ) \in B \st  4\cdot 2^{-n} - 2\cdot 4^{-n} \le r \le
4\cdot 2^{-n} + 2\cdot4^{-n}\}. $$
Let $L_n$ be the tube of leaves of ${\cal V}$ with endpoints in $E_n$. The
intersection of each $L_n$ with the $\theta = 1$ section is also shown in
Figure~\ref{f1dim}, along with the regions $\sigma(E_n)$.  The
figure suggests, and a computation shows, that
$\sigma^{-1}(L_n) \subset \bigcup_k E_k$ for every $n$.  Let
$E_{0,n} = E_n$ and let $L_{0,n} = L_n$.
The set $\bigcup \sigma^{-1}(L_n)$ is also a union of disjoint disks, which we
denote $E_{1,1},E_{1,2}, \ldots\, $.  Let $L_{1,n}$ be the tube of leaves of
${\cal V}$ which begins at $E_{1,n}$.  In general, for each $k$, let $E_{k+1,1},
E_{k+1,2}, \ldots$ be a sequence of disjoint closed disks whose union is
$\sigma^{-1}(\bigcup_n L_{k,n})$, and let $L_{k+1,n}$ be the tube of leaves
beginning at $E_{k+1,n}$. Another computation shows that the diameter of
$E_{k,n}$ goes to zero as $k$ goes to infinity, irrespective of the behavior of
$n$ as a function of $k$.  Moreover, Figure~\ref{f1dim} shows that there are
infinitely many (and in particular more than one) disks $E_{k+1,n'}$ in a given
$E_{k,n}$.  Therefore the closure of the intersection
$$\bigcap_k \bigcup_n E_{k,n}$$
is a Cantor set $C$.

\f1dim

In the leaf $\tilde{T}$, the leaf $T$ is interrupted by the leaf $l$ beginning
at $(0,\frac13,-3)$.  Since the leaf $l$ is infinite, the leaf $\tilde{T}$ never
returns to $T$.  The leaf $l$, in turn, is interrupted by an infinite sequence
of internal entry points, and the $n$th such point lies in $\sigma(D_n)$.  It
follows that the closure of all internal entry points, which is a cross-section
of the closure of $\tilde{T}$, is the Cantor set $C$.  In particular, the
closure of $\tilde{T}$ is 1-dimensional.

The main extra step in the full construction is to take the double-cover $\cal
V'$ of ${\cal V}$.  Similarly, let $V'_s$, $\bar{V'}_s$, and $F'$ be double-covers
of $V_s$, $\bar{V}_s$, and $F$.  Let $\pi:V' \to V$ be the covering map, let
$\beta:V' \to V'$ be the non-trivial deck translation (also use $\pi$ and
$\beta$ for the covering  $F'$ of $F$), and let $\gamma:V' \to V'$ be the
mirror-image involution which switches  $V'_s$ and $\bar{V'}_s$.  Let
$\sigma_1:D_1 \to V'$ be a lift of $\sigma:D \to V'$, where the disk $D_1$ is a
lift of the disk $D$, let $D_2 = \beta(D_1)$, and let $\sigma_2:D_2 \to V'$ be
given by
$$\sigma_2 = \beta \circ \gamma \circ \sigma_1 \circ \beta.$$
Let $D' = D_1 \cup D_2$, and define the self-insertion map $\sigma':D' \to V'$
to be both $\sigma_1$ and $\sigma_2$. The entire system of disks $\{E_{k,n}\}$
and tubes $\{L_{k,n}\}$ has a lift $E_{1,k,n}$ and $L_{1,k,n}$ and another lift
$E_{2,k,n}$ and $L_{2,k,n}$.  Moreover, each $L_{i,k,n}$ intersects
$\im(\sigma')$ in some collection of disks $\{E_{j,k+1,m}\}$, where $j$ and $m$
may both vary.  So the closure of $\tilde{T}$ still has a Cantor set
cross-section and is still 1-dimensional, but since the self-insertion at
$\sigma'$ breaks both circles $T$ and $\bar{T}$, the closure of $\tilde{T}$ is
in this case the minimal set.
\end{proof}

The construction of Theorem~\ref{th1dplug} gives rise to interesting symbolic
dynamics.  Consider the sequence of pairs of integers $(j,n)$ such that the
$i$th internal entry of the leaf $\tilde{T}$ after the interruption of $T$ is in
the disk $E_{j,0,n}$.  The sequence runs:
$$(1,2),(1,1),(1,1),(2,1),(2,1),(1,4),(1,2),(1,1),(1,1),(2,1),(2,1),(1,3),$$
$$(1,2),(1,1),(1,1),(2,1),(2,1),(1,2),(1,1),(1,1),(2,1),(2,1),(1,1),\ldots$$
The sequence is generated by calling {\sc followdisks}(1,$\infty$), where
the recursive procedure {\sc followdisks} is defined as:

\begin{itemize}
\item Procedure {\sc followdisks}($j$,$n$)
\item[1.] Print $(j,n)$ if $n$ is finite.
\item[2.] Do the following three steps with $i = 1,2$:
\item[3.] For each $k$ from $1$ to $n-2$, do {\sc followdisks}($i$,$k$)
if $j+k$ is odd.
\item[4.] Do {\sc followdisks}($i$,$n-1$) twice.
\item[5.] For each $k$ from $n-2$ to $1$, do {\sc followdisks}($i$,$k$)
if $j+k$ is even.
\item[6.] Quit.
\end{itemize}

Finally, note that there are many different constructions that resemble the
construction of Theorem~\ref{th1dplug}, and in general they have similar but
different symbolic dynamics.  For example, where {\sc followdisks} says ``Do
{\sc followdisks}($i$,$n-1$) twice,'' its analogue for some other self-inserted
plug may call itself once or three times or any other integer.

The above PL construction yielding a 1-dimensional minimal set easily generalizes to
$n$-dimensional manifolds, $n\geq 3$: for any $k=1,2,\ldots ,n-2$, there is an aperiodic
PL plug with one minimal set which is of dimension $k$.

\end{document}